\documentclass[fleqn]{mat01}
\usepackage{times,mathtimy,amssymb,latexsym,boxit}
\begin{document}

\setcounter{page}{321}
\firstpage{321}

\font\xx=msam5 at 10pt
\def\ab{\mbox{\xx{\char'03}}}

\newcommand{\R}{\mathbb{R}}
\newcommand{\T}{\mathbb{T}}
\newcommand{\N}{\mathbb{N}}
\newcommand{\ZZ}{\mathbb{Z}}
\newcommand{\C}{\mathbb{C}}

\newtheorem{theo}{Theorem}
\renewcommand\thetheo{\arabic{section}.\arabic{theo}}
\newtheorem{theor}[theo]{\bf Theorem}
\newtheorem{lem}[theo]{Lemma}
\newtheorem{propo}[theo]{\rm PROPOSITION}
\newtheorem{rmk}[theo]{Remark}
\newtheorem{defn}[theo]{\rm DEFINITION}
\newtheorem{exam}{Example}
\newtheorem{coro}[theo]{\rm COROLLARY}
\newtheorem{pot}[theo]{\it Proof of Theorem}

\renewcommand{\theequation}{\thesection\arabic{equation}}

\title{Probabilistic representations of solutions to the heat equation}

\markboth{B Rajeev and S Thangavelu}{Probabilistic representation of the heat equation}

\author{B RAJEEV and S THANGAVELU}

\address{Indian Statistical Institute,
R.V. College Post, Bangalore~560~059, India\\
\noindent E-mail: brajeev@isibang.ac.in; veluma@isibang.ac.in}

\volume{113}

\mon{August}

\parts{3}

\Date{MS received 16 November 2002}

\begin{abstract}
In this paper we provide a new (probabilistic)
proof of a classical  result in partial differential equations,
viz. if $\phi$ is a tempered distribution, then the solution
of the heat equation for the Laplacian, with initial condition
$\phi$, is given by the convolution of $\phi$ with the heat
kernel (Gaussian density). Our results also extend the probabilistic
representation of solutions of the heat equation to initial
conditions that are arbitrary tempered distributions.
\end{abstract}


\maketitle

\section{Introduction}

Let $(X_t)_{t \geq 0}$ be a $d$-dimensional Brownian motion, with $X_0
\equiv 0$. Let $\varphi \in {\cal S}^\prime(\R^d)$, the space of
tempered distributions. Let $\varphi_t$ represent the unique solution to
the heat equation with initial value $\varphi$, viz.
\begin{equation*}
\partial_t \varphi_t = \frac{1}{2} \Delta \varphi_t ~~0 \leq t \leq T;
~~~~~\varphi_0 = \varphi.
\end{equation*}
It is well-known that $\varphi_t =\varphi \ast p_t$, where
$p_t(x) = \frac{1}{(2\pi t)^{d/2}} {\rm e}^{-(|x|^{2}/2t)}$ and `$\ast$'
denotes convolution. When $\varphi$ is smooth, say $\varphi \in {\cal
S}$, the space of rapidly decreasing smooth functions, then the
probabilisitc representation of the solution is given by the equality
$\varphi(t,x) = E \varphi (X_t +x)$ and is obtained by taking
expectations in the Ito formula
\begin{equation*}
\varphi(X_t +x) = \varphi(x) + \int_0^{t} \nabla \varphi (X_s
+x) \cdot {\rm d}X_s + \frac{1}{2} \int_0^{t} \Delta \varphi
(X_s+x) {\rm d}s.
\end{equation*}

Such representations are well-known (see \cite{1,2,3,4}) and extend to a large
class of initial value problems, with the Laplacian $\Delta$ replaced by a
suitable (elliptic) differential operator $L$ and $(X_t)$ being replaced
by the diffusion generated by $L$. A basic problem here is to extend the
representation to situations where $\varphi$ is not smooth.

The main contribution of this paper is to give a probabilistic
representation of solutions to the initial value problem for the
Laplacian with an arbitrary initial value $\varphi \in {\cal S}^\prime$.
This representation follows from the Ito formula developed in \cite{9}, for
the ${\cal S}^\prime$-valued process $(\tau_{X_t} \varphi)$, where
$\tau_x \varphi$ is the translation of $\varphi$ by $x \in \R^d$. Our
representation (Theorem~2.4) then reads, $\varphi_t =E \tau_{X_t}
\varphi$ where of course $\varphi_t$ is the solution of the initial
value problem for the Laplacian, with initial value $\varphi \in {\cal
S}^\prime$. In particular, the fundamental solution $p_t (x- \cdot)$ has
the representation, $p_t (x-\cdot) = E \tau_{X_t} \delta_x$. However, the
results of \cite{9} only show that if $\varphi \in {\cal S}^\prime_p$, then
there exists $q > p$ such that the process $(\tau_{X_t} \varphi)$ takes
values in ${\cal S}^\prime_q$. Here for each real $p$, the ${\cal
S}_p$s are the `Sobolev spaces' associated with the spectral
decomposition of the operator $|x|^2-\Delta$ or equivalently they are
the Hilbert spaces defining the countable Hilbertian structure of ${\cal
S}^\prime$ (see~\cite{6}). ${\cal S}^\prime_p$, the dual of ${\cal S}_p$, is
the same as ${\cal S}_{-p}$. Clearly it would be desirable to have the
process $(\tau_{X_t}\varphi)$ take values in ${\cal S}^\prime_p$,
whenever $\varphi \in {\cal S}^\prime_p$. Such a result also has
implications for the semi-martingale structure of the process
$(\tau_{X_t})$ -- it is a semi-martingale in ${\cal S}^\prime_{p+1}$
(Corollary~2.2) and fails to have this property in ${\cal S}^\prime_q$
for $q < p+1$ (see Remark 5.2 of \cite{5}).

Given the above remarks and the results of \cite{9}, the properties of the
translation operators become significant. We show in Theorem~2.1 that
the operators $\tau_x : {\cal S}_p \rightarrow {\cal S}_p$ for $x \in
\R^d$, are indeed bounded operators, for any real $p$, with the operator
norms being bounded above by a polynomial in $|x|$. The proof uses
interpolation techniques well-known to analysts. Theorem~2.4 then
gives a comprehensive treatment of the initial value problem for the
Laplacian from a probabilistic point of view.

\section{Statements of the main results}

\setcounter{equation}{0}

Let $(\Omega, {\cal F}, ({\cal F}_t)_{t \geq 0}, P)$ be a filtered probability
space with a filtration $({\cal F}_t)$ satisfying usual conditions:
${\cal F}_t = \bigcap_{s > t} {\cal F}_s$ and ${\cal F}_0$ contains
all $P$-null sets. Let $(X_t)_{t \geq 0}$ be a $d$-dimensional,
$({\cal F}_t)$-Brownian motion with $X_0 \equiv 0$.

${\cal S}$ denotes the space of rapidly decreasing smooth functions on $\R^d$
(real valued) and ${\cal S}^\prime$ its dual, the space of tempered distributions.
We refer to \cite{11} for formal definitions. For $x \in \R^d, \delta_x \in
{\cal S}^\prime$ will denote the Dirac distribution at $x$.
Let $\{\tau_x : x \in \R^{d}\}$ denote the translation operators defined
on functions by the formula $\tau_x f(y)=f(y-x)$ and let $\tau_x :
{\cal S}^\prime \rightarrow {\cal S}^\prime$ act on distributions by
\begin{equation*}
\langle \tau_x \varphi, f \rangle = \langle \varphi, \tau_{-x} f\rangle.
\end{equation*}
The nuclear space structure of ${\cal S}^\prime$ is given by the family of Hilbert
spaces ${\cal S}_p, p \in \R$, obtained as the completion of
${\cal S}$ under the Hilbertian
norms $\{\| \cdot\|_p\}_{p \in \R}$ defined by
\begin{equation*}
\|\varphi\|_p^2 =\sum\limits_k (2|k| +d)^{2p} \langle \varphi, h_k \rangle^2,
\end{equation*}
where $\varphi \in {\cal S}$, and the sum is taken over $k=(k_1, \ldots, k_d) \in
\ZZ^d_+, |k|=(k_1 + \cdots + k_d), \langle \varphi, h_k \rangle$ denotes the
inner product in $L^2(\R^d)$ and $\{h_k, k \in \ZZ^d_+\}$ is the ONB
in $L^2(\R^d)$, constructed as follows: for $x=(x_1, \ldots, x_d), h_k(x)
=h_{k_1} (x_1) \ldots h_{k_d} (x_d)$.
The one-dimensional Hermite functions are given by $h_\ell (s)=\frac{1}
{(\sqrt \pi 2^\ell \ell!)^{1/2}} {\rm e}^{-(s^2/2)} H_\ell (s)$, where
$H_\ell (s) =(-1)^\ell {\rm e}^{s^2} \frac{{\rm d}^\ell}{{\rm d}s^\ell} {\rm e}^{-s^2}$ are the
Hermite polynomials. While we mainly deal with real valued functions, at times we need to use
complex valued functions. In such cases, the spaces ${\cal S}_p$ are defined in
a similar fashion as above, i.e. as the completion of ${\cal S}$ with respect to $\|
\cdot \|_p$. However, in the definition of $\|\varphi\|^2_p$ above we need
to replace the real $L^2$ inner product $\langle \varphi, h_k \rangle$ by the one
for complex valued functions, viz. $\langle \varphi, \psi\rangle = \int_{\R^d}
\varphi (x) \bar \psi (x) {\rm d}x$ and $\langle \varphi, h_k \rangle^2$ is replaced by
$| \langle \varphi, h_k \rangle|^2$. It is well-known (see \cite{6,7}) that ${\cal S}
=\bigcap_ p {\cal S}_p, {\cal S}^\prime =\bigcup_p {\cal S}_p$
and ${\cal S}^\prime_p=:$ dual of ${\cal S}_p={\cal S}_{-p}$.
We will denote by $\langle \cdot, \!\cdot \rangle_p$, the inner product corresponding to
the norm $\| \cdot\|_p$.

Let $(Y_t)_{t \geq 0}$ be an ${\cal S}_p$-valued, locally bounded, previsible
process, for some $p \in \R$. Let $\partial_i : {\cal S}_p
\rightarrow {\cal S}_{p-1/2}$
be the partial derivatives, $1 \leq i \leq d$, in the sense of distributions.
Then since $\partial_i, 1 \leq i \leq d$ are bounded linear operators it
follows that $(\partial_i Y_t)_{t \geq 0}$ is an ${\cal S}_{p-1/2}$-valued,
locally bounded, previsible process. From the theory of stochastic integration
in Hilbert spaces \cite{8}, it follows that the processes
\begin{equation*}
\left( \int_0^{t} Y_s {\rm d} X^i_s \right)_{t \geq 0}, \left(
\int_0^{t} \partial_i Y_s {\rm d}X_s^i \right)_{t \geq 0}
\end{equation*}
are continuous ${\cal F}_t$ local martingales for $1 \leq i \leq d$, with
values in ${\cal S}_p$ and ${\cal S}_{p-1/2}$ respectively. If $X_t =(X^1_t, \ldots,
X_t^d)$ is a continuous $\R^d$-valued, ${\cal F}_t$-semi-martingale,
it follows from the general theory that the above processes too are
continuous ${\cal F}_t$-semi-martingales with values in
${\cal S}_p$ and ${\cal S}_{p-1/2}$
respectively.

\begin{theor}[\!]
Let $p \in \R$. There exists a polynomial $P_k(\cdot)$ of degree
$k=2([|p|]+1)$ such that the following holds{\rm :} For $x \in \R^{d}, \tau_x
:{\cal S}_p \rightarrow {\cal S}_p$ is a bounded linear map and we have
\begin{equation*}
\|\tau_x \varphi\|_p \leq P_k (|x|) \|\varphi\|_p
\end{equation*}
for all $\varphi \in {\cal S}_p$.
\end{theor}

In (\cite{9}, Theorem~2.3) we showed that if $(X_t)_{t\geq 0}$ is a continuous,
$d$-dimensional, ${\cal F}_t$-semi-martingale and $\varphi \in {\cal S}_p
\subset {\cal S}^\prime$, then the process $(\tau_{X_t} \varphi)_{t \geq 0
}$ is an ${\cal S}_q$-valued continuous semi-martingale for some $q < p$.
Corollary~2.2 below says that we can take $q=p-1$.

\begin{coro}$\left.\right.$\vspace{.5pc}

\noindent Let $(X_t)_{t \geq 0}$ be a continuous $d$-dimensional{\rm ,}
${\cal F}_t$-semi-martingale. Let $\varphi \in {\cal S}_p, p \in \R$. Then $(\tau_{X_t}
\varphi)_{t \geq 0}$ is an ${\cal S}_p$-valued{\rm ,} continuous adapted process.
Moreover it is an ${\cal S}_{p-1}$-valued{\rm ,} continuous ${\cal F}_t$-semi-martingale
and the following Ito formula holds in ${\cal S}_{p-1}${\rm :} a.s.{\rm ,} $\forall~ t \geq 0${\rm ,}
\begin{align}
\tau_{X_t} \varphi &= \tau_{X_0} \varphi - \sum\limits^d_{i=1} \int_0^{t} \partial_i (\tau_{X_s} \varphi) {\rm d} X_s^i \nonumber \\
&\quad\ + \frac{1}{2} \sum\limits^d_{i,j=1} \int_0^{t} \partial^2_{ij} (\tau_{X_s} \varphi) {\rm d} \langle X^i,
X^j \rangle_s,
\end{align}
where $X_t =(X^1_t, \ldots, X_t^d)$ and $(\langle X^i, X^j \rangle_t)$ is the quadratic
variation process between $(X^i_t)$ and $(X^j_t), 1 \leq i, j \leq d$.
\end{coro}

\begin{proof}
From Theorem~2.1, it follows that $(\tau_{X_t} \varphi)$ is an
${\cal S}_p$-valued continuous adapted process. By Theorem 2.3 of \cite{9}, $\exists~
q < p$, such that $(\tau_{X_t} \varphi)$ is an ${\cal S}_q$ semi-martingale
and the above equation holds in ${\cal S}_q$. Clearly each of the terms in the
above equation is in ${\cal S}_{p-1}$ and the result follows.\hfill \ab
\end{proof}

The next corollary pertains to the case when $(X_t) =(X_t^1, \ldots, X_t^d)$
is a $d$-dimensional Brownian motion, $X_0 \equiv 0$. In (\cite{5},
Definition~3.1), we introduced the notion of an ${\cal S}^\prime_p (={\cal S}_{-p},
p > 0)$-valued strong solution of the SDE
\begin{align}
&{\rm d}Y_t = \frac{1}{2} \Delta (Y_t) {\rm d}t + \nabla Y_t \cdot {\rm d}X_t, \nonumber\\
&Y_0 = \varphi,
\end{align}
where $\nabla = (\partial_1, \ldots, \partial_d)$ and $\Delta
= \sum_{i=1}^d \partial_i^2$. There we showed that if $\varphi \in {\cal
S}^\prime_p$, then the above equation has a unique
${\cal S}^\prime_q$-valued strong solution, $q \geq p+2$. Theorem~2.1 implies that we
indeed have an (unique) ${\cal S}^\prime_p$-valued strong solution.

\begin{coro}$\left.\right.$\vspace{.5pc}

\noindent Let $\varphi \in {\cal S}^\prime_p$. Then{\rm ,} eq.~{\rm
(2.2)} has a unique ${\cal S}^\prime_p$-valued strong solution on $0
\leq t \leq T$.
\end{coro}

\begin{proof}
By Corollary~2.2, the process $(\tau_{X_t} \varphi)$, where $(X_t)$ is a $d$-dimensional
Brownian motion, $X_0 \equiv 0$, satisfies eq.~(2.1). Further,
\begin{equation*}
E \int_0^{T} \|\tau_{X_t} \varphi\|_{-p}^2 {\rm d}t = \int_0^{T}
\int_{\R^d} \|\tau_x \varphi\|_{-p}^2 \frac{{\rm e}^{-(|x|^2/2t)}}
{(2\pi t)^{d/2}} {\rm d}x~{\rm d}t < \infty.
\end{equation*}
Uniqueness follows as in Theorem~3.3 of \cite{5}.\hfill \ab
\end{proof}

We now consider the heat equation for the Laplacian with initial condition
$\varphi \in {\cal S}_p$, for some $p \in \R$.
\begin{align}
&\partial_t \varphi_t = \frac{1}{2} \Delta \varphi_t ~~0 < t \leq T, \nonumber\\
&\varphi_0 = \varphi.
\end{align}

By an ${\cal S}_p$-valued solution of (2.3), we mean a continuous map $t \rightarrow
\varphi_t\!: [0, T] \rightarrow S_p$ such that the following
equation holds in ${\cal S}_{p-1}$:
\begin{equation}
\varphi_t =\varphi + \int_0^{t} \frac{1}{2} \Delta \varphi_s {\rm d}s.
\end{equation}
Let $\{h_k^{p-1}\}$ be the ONB in ${\cal S}_{p-1}$ given by $h_k^{p-1}=
(2|k|+d)^{-(p-1)}h_k$. We then have for $p<0$ and $t \leq T$:
\begin{align*}
\|\varphi_t\|^2_{p-1} &= \sum\limits^\infty_{|k|=0} \langle \varphi_t,
h_k^{p-1} \rangle^2_{p-1} \\
&= \sum\limits_{|k|=0}^\infty \left\{ \langle \varphi, h_k^{p-1}\rangle^2_
{p-1} + 2 \int_0^{t} \langle \varphi_s, h_k^{p-1} \rangle_{p-1}
{\rm d}\langle \varphi_s, h_k^{p-1} \rangle_{p-1}\right\} \\
&= \|\varphi\|^2_{p-1} + \sum\limits_{|k|=0}^\infty 2 \int_0^{t}
\langle \varphi_s, h_k^{p-1} \rangle_{p-1} \left\langle \frac{1}{2} \Delta
\varphi_s, h_k^{p-1}\right\rangle_{p-1} {\rm d}s\\
&= \|\varphi\|^2_{p-1} + 2 \int_0^{t} \left\langle \frac{1}{2}
\Delta \varphi_s, \varphi_s \right\rangle_{p-1} {\rm d}s.
\end{align*}
It follows from the results of \cite{5} (the monotonicity condition) that
for $p <0$,
\begin{equation*}
2 \left\langle \frac{1}{2} \Delta \varphi, \varphi \right\rangle_{p-1} +\sum\limits_{i=1}
^d \|\partial_i \varphi\|^2_{p-1} \leq C~\|\varphi\|^2_{p-1}
\end{equation*}
for some constant $C >0$ for all $\varphi \in {\cal S}_p$. We then get
\begin{equation*}
\|\varphi_t\|^2_{p-1} \leq \|\varphi\|^2_{p-1} + C\int_0^{t}
\|\varphi_s\|^2_{p-1} {\rm d}s.
\end{equation*}
Hence for the case $p <0$, uniqueness
follows from the Gronwall lemma. Uniqueness for the case $p \geq 0$,
follows from uniqueness for the case $p < 0$ and the inclusion
${\cal S}_p \subset {\cal S}_q$ for $q < p$. It is well-known that the
solutions of the initial value problem~(2.3) in ${\cal S}^\prime(\R^d)$ are
given by convolution of $\varphi$
and $p_t(x)$, the heat kernel. That these  coincide (as they should) with
the ${\cal S}_p$-valued solutions follows from the `probabilistic representation'
given by Theorem~2.4 below. Define the Brownian semi-group $(T_t)_{t \geq 0}$
on ${\cal S}$ in the usual manner:
\begin{equation*}
T_t \varphi (x) = \varphi \ast p_t (x) ~~t > 0,~~~~~~ T_0 \varphi = \varphi
\end{equation*}
where $p_t (x) = \frac{1}{(2\pi t)^{d/2}} {\rm e}^{(-|x|^2/2t)},
t > 0$ and `$\ast$' denotes convolution: $f \ast g (x) =\int_{\R^d} f(y)
g(x-y) {\rm d}y$. In the next theorem we consider standard Brownian motion $(X_t)$.

\begin{theor}[\!]
{\rm (a)} Let $\varphi \in {\cal S}_p$. Then for $t \geq 0${\rm ,} the ${\cal S}_p$-valued
random variable $\tau_{X_t} \varphi$ is Bochner integrable and we have
\begin{equation*}
E ~\tau_{X_t} \varphi = \varphi \ast p_t  = T_t \varphi.
\end{equation*}
In particular{\rm ,} for every $p \in \R${\rm ,} and $T > 0, \sup_{t \leq T}
\|T_t\| < \infty$ where $\|T_t\|$ is the operator norm of $T_t: {\cal S}_p
\rightarrow {\cal S}_p$.\\
{\rm (b)} For $\varphi \in {\cal S}_p${\rm ,} the initial value problem~{\rm (2.3)} has
a unique ${\cal S}_p$-valued solution $\varphi_t$ given by
\begin{equation*}
\varphi_t = E \tau_{X_t} \varphi.
\end{equation*}
Further $\varphi_t \rightarrow \varphi$ strongly in ${\cal S}_p$ as
$t \rightarrow 0$.
\end{theor}

\section{Proofs of Theorems 2.1 and 2.4}

\setcounter{theo}{0}
\setcounter{equation}{0}
The spaces ${\cal S}_p$ can be described in terms of the spectral
properties of
the operator $H$ defined as follows:
\begin{equation*}
Hf = (|x|^2 -\Delta) f,~~f \in {\cal S}.
\end{equation*}
If $\{h_k\}$ is the ONB in $L^2 (\R^d)$ consisting of Hermite functions
(defined in \S2), then it is well-known (see \cite{10}) that
\begin{equation*}
H h_k = (2 |k| +d) h_k.
\end{equation*}
For $f \in {\cal S}$, define the operator $H^p$ as follows:
\begin{equation*}
H^p f = \sum\limits_k (2|k|+d)^p \langle f, h_k \rangle h_k.
\end{equation*}
Here $p$ is any real number. For $f \in {\cal S}$ and $z=x+iy \in ~\C$ define
$H^z f = \sum_k (2 |k| +d)^z \langle f, h_k\rangle h_k$ and note that, $H^z f=
H^x (H^{iy} f) = H^{iy} (H^x f)$ and $H^{iy} : L^2 \rightarrow L^2$ is an
isometry. Further,
\begin{align*}
\|H^z f \|_0^2 &= \sum\limits_k (2|k| +d)^{2x} \langle f, h_k\rangle^2 \\
&= \|f\|_x^2.
\end{align*}
The following propositions (3.1, 3.2 and 3.3) may be well-known. We
include the proofs for completeness.

\begin{propo}$\left.\right.$\vspace{.5pc}

\noindent For any $p$ and $q, \|H^p \varphi\|_{q-p} =
\|\varphi\|_q$ for $\varphi \in {\cal S}$. Consequently{\rm ,} $H^p :
{\cal S}_q \rightarrow {\cal S}_{q-p}$ extends as a linear isometry.
Moreover, this isometry is onto.
\end{propo}

\begin{proof}
Let $h_k^p = (2|k| +d)^{-p} h_k$. Then from the relation $\langle
\varphi, h_k \rangle_p = (2|k|+d)^{2p} \langle\varphi,h_k\rangle$ it follows that $\{h_k^p\}$
is an ONB for ${\cal S}_p$. Let $\varphi \in {\cal S}$. Since
\begin{align*}
H^p \varphi &= \sum\limits_k \langle\varphi, h_k\rangle (2|k|+d)^p h_k \\
&= \sum\limits_k \langle\varphi, h_k\rangle (2|k| +d)^q h_k^{q-p} ,
\end{align*}
we get $\|H^p \varphi\|^2_{q-p} = \|\varphi\|^2_q$.

To show that $H^p$ is onto, consider $\psi \in {\cal S}_{q-p}$,
\begin{equation*}
\psi = \sum_k \langle \psi, h_k^{q-p} \rangle_{q-p} h_k^{q-p}.
\end{equation*}
Defining $\varphi =: \sum_k \langle\psi, h_k^{q-p}\rangle_{q-p}
h_k^q,$ we see that $\varphi \in {\cal S}_q$. Also,\vspace{.2pc}
\begin{equation*}
H^p \varphi = \sum\limits_k \langle\varphi, h_k^q\rangle_q h_k^{q-p}
= \sum\limits_k \langle\psi, h_k^{q-p}\rangle_{q-p} h_k^{q-p}
= \psi.
\end{equation*}

$\left.\right.$\vspace{-2.7pc}

\hfill \ab\vspace{1pc}
\end{proof}

Let $A_j = x_j +\partial_j$ and $A_j^+ = x_j - \partial_j,
1 \leq j \leq d$.
Then it is easy to see that
\begin{equation*}
H = \frac{1}{2} \sum\limits_{j=1}^d (A_j A_j^+ + A_j^+ A_j).
\end{equation*}
For multi-indices $\alpha = (\alpha_1, \ldots, \alpha_d)$ and $\beta = (\beta_1,
\ldots, \beta_d)$ we define
\begin{equation*}
A^\alpha =: A_1^{\alpha_1} \ldots A_d^{\alpha_d},~~~~~~~ (A^+)^\beta
=: (A_1^+)^{\beta_1} \ldots (A_d^+)^{\beta_d}.
\end{equation*}
For an integer $\ell \geq 0$ and $x \in \R$, recall that
\begin{equation*}
h_\ell (x) = \frac{1}{(\sqrt \pi 2^\ell \ell!)^{1/2}} {\rm e}^{-(x^2/2)}
H_\ell (x),
\end{equation*}
where $H_\ell$ is the Hermite polynomial defined by
\begin{equation*}
H_\ell (x) = (-1)^\ell {\rm e}^{x^2} \frac{{\rm d}^\ell}{{\rm d}x^\ell}
{\rm e}^{-x^2}.
\end{equation*}
It is easily verified that
\begin{align*}
\left(x+\frac{\rm d}{{\rm d}x} \right) \left( {\rm e}^{-(x^2/2)} H_\ell (x) \right)
&= 2\ell \left( {\rm e}^{-(x^2/2)} H_{\ell-1} (x) \right), \\
\left( x-\frac{\rm d}{{\rm d}x} \right) \left( {\rm e}^{-(x^2/2)} H_\ell (x) \right) &=
{\rm e}^{-(x^2/2)} H_{\ell+1} (x).
\end{align*}
It then follows that
\begin{align*}
A_j^+ h_{k_j} (x_j) &= \sqrt {2(k_j +1)} h_{{k_j}+1} (x_j), \\
A_j h_{k_j} (x_j) &= \sqrt {2 k_j} h_{{k_j}-1} (x_j).
\end{align*}
Iterating these two formulas we get the following:

\begin{propo}$\left.\right.$\vspace{.5pc}

\noindent Let $k, \beta$ and $\alpha$ be multi-indices such that $k_j \geq \alpha_j,
j=1, \ldots, d$. Then
\begin{align*}
&(A^+)^\beta h_k (x) = 2^{|\beta|/2} \left( \frac{(k+\beta)!}{k!} \right)^{1/2} h_{k+\beta} (x),\\
&A^\alpha h_k (x) = 2^{|\alpha|/2} \left( \frac{k!}{(k-\alpha)!}
\right)^{1/2} h_{k-\alpha}(x),
\end{align*}
where $k! = k_1 ! \ldots k_d!$.
\end{propo}

\begin{propo}$\left.\right.$\vspace{.5pc}

\noindent For all $m \geq 0, \exists$ constants $C_1 =C_1 (m)$ and
$C_2 =C_2(m)$ such that the following hold{\rm :}\\
{\rm (a)} For all $f \in {\cal S}${\rm ,}
\begin{equation*}
\|f\|_m \leq C_1~\sum\limits_{|\alpha| +|\beta| \leq 2m} \|A^\alpha
(A^+)^\beta f\|_0 \leq  C_2 \|f\|_m.
\end{equation*}
{\rm (b)} For all $f \in {\cal S}${\rm ,}
\begin{equation*}
\|f\|_m \leq C_1 \sum\limits_{|\alpha|+|\beta|\leq 2m} \|x^\alpha \partial^
\beta f\|_0 \leq C_2 \|f\|_m.
\end{equation*}
\end{propo}

\begin{proof}
(a) We can write
\begin{equation*}
H^m = \sum\limits_{|\alpha| +|\beta|\leq 2m} C_{\alpha \beta} A^\alpha
(A^+)^\beta,
\end{equation*}
where $C_{\alpha \beta}$ are constants. Since $\|f\|_m = \|H^m f\|_0$,
the first part of the inequality follows.
To show the second half of the inequality it is sufficient to show that
for $f \in {\cal S}$ and $|\alpha| +|\beta| \leq 2m, \|A^\alpha (A^+)^\beta
H^{-m} f \|_0 \leq C_{\alpha \beta} \|f\|_0$. Now,
\begin{align*}
\|A^\alpha (A^+)^\beta H^{-m} f\|_0^2 &= \sum\limits_\ell
\langle A^\alpha (A^+)^\beta H^{-m} f, h_\ell \rangle^2 \\
&= \sum\limits_\ell \left[ \sum\limits_k (2|k|+d)^{-m} \langle f, h_k\rangle \langle A^\alpha
(A^+)^\beta h_k, h_\ell \rangle\right]^2 \\
&= \sum\limits_\ell \left[\sum\limits_k (2|k| +d)^{-m} \langle f, h_k \rangle
C_{k,\beta,\alpha} \langle h_{k+\beta-\alpha}, h_\ell \rangle \right]^2 \\
&= \sum\limits_\ell (2|\ell +\alpha\!-\!\beta|+d)^{-2m} C^2_{\ell +\alpha-\beta,
\beta,\alpha} \langle f, h_{\ell+\alpha-\beta} \rangle^2,
\end{align*}
where the sum is taken over $\ell=(\ell_1, \ldots, \ell_d)$ such that
$\ell_j +\alpha_j -\beta_j \geq 0$ for $1 \leq j \leq d$ and where we have
used Proposition~3.2 in the last but one equality above. From the same proposition,
it follows that
\begin{equation*}
(2 |\alpha +\ell -\beta| +d)^{-2m} C^2_{\ell+\alpha-\beta, \beta,\alpha}
\end{equation*}
are uniformly bounded in $\ell$ for $|\alpha| +|\beta| \leq 2m$ and the
second inequality in (a) follows.

\noindent (b) Since $\|f\|_m =\|H^m f\|_0$ and clearly $H^m =\sum_{|\alpha|+
|\beta| \leq 2m} C_{\alpha \beta} x^\alpha \partial^\beta$, the first
inequality follows. To prove the second inequality, note that
\begin{equation*}
x_j = \frac{1}{2} (A_j + A_j^+),~~~~~~ \partial_j = \frac{1}{2} (A_j - A_j^+).
\end{equation*}
Hence, using $[A_j, A^+_k]=\delta_{jk}I$,
\begin{equation*}
x^\alpha \partial^\beta = \sum\limits_{|k|+|\ell|\leq |\alpha|+|\beta|}
C_{k,\ell} A^k (A^+)^\ell
\end{equation*}
and hence by part (a) we get
\begin{equation*}
\sum\limits_{|\alpha| +|\beta| \leq 2m} \|x^\alpha \partial^\beta f\|_0
\leq C_1 \sum\limits_{|k| +|\ell| \leq 2m} \|A^k (A^+)^\ell f\|_0
\leq C_2 \|H^m f\|_0.
\end{equation*}
\hfill \ab
\end{proof}

\setcounter{section}{2}
\setcounter{theo}{0}
\begin{pot}
{\rm We first show that for an integer $m \geq 0$,
\begin{equation*}
\|\tau_x \varphi\|_m \leq P_{2m} (|x|) \|\varphi\|_m,
\end{equation*}
where $P_{2m} (t)$ is a polynomial in $t \in \R$ of degree $2m$ with non-negative
coefficients. This follows from Proposition~3.3:
\begin{align*}
\|\tau_x f\|_m &\leq C_1 \sum\limits_{|\alpha| +|\beta| \leq 2m} \|y^\alpha \partial^\beta
\tau_x f\|_0 \\
&\leq C_1 \sum\limits_{|\alpha|+|\beta| \leq 2m} \|(y+x)^\alpha \partial^
\beta f\|_0.
\end{align*}
The last sum is clearly dominated by $P_{2m} (|x|) \|f\|_m$ for some polynomial
$P_{2m}$. If $m < p < m+1$, where $m \geq 0$ is an integer, we prove the result using the
3-line lemma: for $f, g \in {\cal S}$, let
\begin{equation*}
F(z) = \langle H^{\bar z} \tau_x H^{-z} f, g \rangle_0.
\end{equation*}
Then from the expansion in $L^2$ for the RHS it is verified that $F(z)$ is
analytic in $m < {\rm Re}~ z < m + 1$ and continuous in $m \leq {\rm
Re}~ z \leq m+1$. We will show that
\setcounter{section}{3}
\begin{align}
&|F(m+iy)| \leq P_{2m} (|x|) \|f\|_0 \|g\|_0,\nonumber\\[.2pc]
&|F(m+1+iy)| \leq P_{2(m+1)} (|x|) \|f\|_0 \|g\|_0
\end{align}
for $-\infty < y < \infty$. Hence from the 3-line lemma \cite{12}, it follows
that
\begin{align*}
|F(p+iy)| &\leq (P_{2m} (|x|) \|f\|_0 \|g\|_0)^{m+1-p} (P_{2(m+1)}
(|x|) \|f\|_0 \|g\|_0)^{p-m}  \\[.2pc]
&\leq P_k (|x|) \|f\|_0 \|g\|_0,
\end{align*}
where $P_k (t)$ is a polynomial in $t$ of degree $k=2([p]+1)$. It follows
that
\begin{equation*}
\|\tau_x f\|_p \leq P_k (|x|) \|f\|_p.
\end{equation*}
Using the fact that ${\cal S}_{-p} ={\cal S}^\prime_p$ we get $\|\tau_x f\|_{-p} \leq P_k
(|x|) \|f\|_{-p}$ for $m \leq p \leq m+1$.

The following chain of inequalities establish the inequalities (3.1):
\begin{align*}
|F(m+iy)| &\leq \|H^{m-iy} \tau_x H^{-(m+iy)}f\|_0 \|g\|_0\\[.2pc]
&\leq \|H^m \tau_x H^{-(m+iy)} f\|_0 \|g\|_0 \\[.2pc]
&= \|\tau_x H^{-(m+iy)} f\|_m \|g\|_0 \\[.2pc]
&\leq P_{2m} (|x|) \|H^{-(m+iy)} f\|_m \|g\|_0 \\[.2pc]
&= P_{2m} (|x|) \|H^{-iy} f\|_0 \|g\|_0 \\[.2pc]
&= P_{2m} (|x|) \|f\|_0 \|g\|_0.
\end{align*}
This completes the proof of Theorem~2.1.}\hfill \ab
\end{pot}

\setcounter{section}{2}
\setcounter{theo}{3}
\begin{pot}
{\rm (a) Let $\varphi \in {\cal S}_p, p \in \R$. From Theorem 2.1 we have
\setcounter{section}{3}
\begin{equation*}
\|\tau_{X_t} \varphi\|_{p} \leq P_k (|X_t|) \|\varphi\|_p,
\end{equation*}
where $P_k$ is a polynomial. Since $EP_k (|X_t|) < \infty $,
Bochner integrability follows. For $\psi \in {\cal S}, \varphi \in {\cal S}$,
\begin{align*}
\left\langle \psi, \int \tau_x \varphi~ p_t (x) {\rm d}x \right\rangle &= \int \langle
\psi, \tau_x \varphi \rangle p_t (x) {\rm d}x \\
&= \int p_t (x) {\rm d}x \int \psi (y) \varphi (y-x) {\rm d}y \\
&= \int \psi (y) {\rm d}y \int \varphi (y-x) p_t (x) {\rm d}x \\
&= \int \psi (y) \varphi \ast p_t (y) {\rm d}y \\
&= \langle \psi, \varphi \ast p_t\rangle.
\end{align*}
The result for $\varphi \in {\cal S}_p$ follows by a continuity argument: Let
$\varphi_n \in {\cal S}, \varphi_n \rightarrow \varphi$ in ${\cal S}_p$. Hence $\varphi_n
\ast p_t \rightarrow \varphi \ast p_t$ weakly in ${\cal S}^\prime$. Hence,
\begin{align*}
\langle \psi, \varphi \ast p_t \rangle &= \lim\limits_{n \rightarrow \infty} \langle \psi,
\varphi_n \ast p_t \rangle\\
&= \lim\limits_{n \rightarrow \infty} \int \psi (y) \varphi_n \ast p_t
(y) {\rm d}y \\
&= \lim\limits_{n \rightarrow \infty} \int \langle \psi, \tau_x \varphi_n\rangle p_t
(x) {\rm d}x \\
&= \int \langle \psi, \tau_x \varphi\rangle p_t (x) {\rm d}x \\
&= \left\langle \psi, \int \tau_x \varphi~ p_t (x) {\rm d}x \right\rangle,
\end{align*}
where we have used DCT in the last but one equality. That $T_t : {\cal S}_p \rightarrow
{\cal S}_p$ is a (uniformly) bounded operator follows:
\pagebreak

$\left.\right.$\vspace{-2pc}

\begin{align*}
\|T_t \varphi\|_p &= \|\varphi \ast p_t\|_p = \|E \tau_{X_t} \varphi\|_p \\
&= \left\|\int \tau_x \varphi~ p_t (x) {\rm d}x \right\|_p
\leq \int \|\tau_x \varphi\|_p p_t (x) {\rm d}x \\
&\leq \|\varphi\|_p \int P_k (|x|) p_t (x) {\rm d}x
\leq C~\|\varphi\|_p,
\end{align*}
where $C = \sup_{s \leq T} \int P_k (|x|) p_s (x) {\rm d}x < \infty$.

\noindent (b) Let  $(X_t)$ be the standard Brownian
motion so that $\langle X^i, X^j \rangle \equiv 0$ for $i \neq j$. Equation~(2.1)
then reads, for $\varphi \in {\cal S}_p, p \in \R$,
\begin{align}
\tau_{X_t} \varphi &= \varphi -\int_0^{t} \nabla (\tau_{X_s}
\varphi) \cdot {\rm d}X_s + \frac{1}{2} \int_0^{t} \Delta (\tau_{X_s} \varphi) {\rm d}s.
\end{align}
The stochastic integral is a martingale in ${\cal S}_{p-1}$:
\begin{align*}
E~ \left\| \int_0^{t} \partial_i (\tau_{X_s} \varphi) {\rm d}X_s^i \right\|^2_{p-1}
&\leq C_1 ~E \int_0^{t} \|\partial_i (\tau_{X_s} \varphi)\|_{p-1}^2
{\rm d}s\\
&= C_1 \int_0^{t} \left( \int \|\partial_i (\tau_x \varphi)\|^2_{p-1} p_s (x) {\rm d}x \right) {\rm d}s \\
&\leq C_2 \int_0^{t} \left( \int \| \tau_x \varphi\|^2_p~ p_s(x) {\rm d}x \right) {\rm d}s \\
&\leq C_3 \|\varphi\|_p \int_0^{t} \left( \int P_k (|x|) p_s (x) {\rm d}x \right) {\rm d}s \\
&< \infty.
\end{align*}
Let $\varphi_t =E \tau_{X_t} \varphi$. Taking expected values in (3.2)
we get eq.~(2.4). Hence $\varphi_t$ is the  solution to the
heat equation with initial value $\varphi \in {\cal S}_p$. The uniqueness
of the solution is well-known and also follows from the remarks preceeding
the statement of Theorem~2.4.

To complete the proof of the theorem, we need to show that $\varphi_t
\rightarrow \varphi$ in ${\cal S}_p$ as $t \downarrow 0$. Let ${\cal F}$ denote
the Fourier transform, i.e. ${\cal F}f (\xi) =\int {\rm e}^{-i (x\cdot \xi)}
f(x) {\rm d}x$ for $f \in {\cal S}$. Then ${\cal F}$ extends to ${\cal S}^\prime$ by duality, where
we consider ${\cal S}^\prime$ as a complex vector space. Since ${\cal F} (h_n)=
(- \sqrt {-1})^n h_n$ (\cite{10}, p.~5, Lemma~1.1.3), ${\cal F}$ acts as a
bounded operator from ${\cal S}_p$ to ${\cal S}_p$, for all $p$. Let $\varphi \in {\cal S}_p$.
\begin{equation*}
\varphi_t -\varphi = T_t \varphi - \varphi
= {\cal F}^{-1} (S_t ({\cal F} \varphi)),
\end{equation*}
where
\begin{equation*}
S_t \varphi (x) = {\cal F} (T_t -I) {\cal F}^{-1} \varphi (x) =
({\rm e}^{-(t/2)|x|^2} -1) \varphi (x).
\end{equation*}

Clearly, $S_t: {\cal S}_p \rightarrow {\cal S}_p$ is a bounded operator and
\begin{equation*}
\|\varphi_t -\varphi\|_p = \|S_t ({\cal F} \varphi)\|_p.
\end{equation*}
The following proposition completes the proof of the theorem.}
\end{pot}

\setcounter{theo}{3}
\begin{propo}$\left.\right.$\vspace{.5pc}

\noindent Let $\varphi \in {\cal S}_p, p \in \R$. Then $\|S_t \varphi\|_p \rightarrow 0$
as $t \rightarrow 0$.
\end{propo}

\begin{proof}
We prove the proposition by showing that (i) $S_t: {\cal S}_p
\rightarrow {\cal S}_p$ are uniformly bounded, $0 < t \leq T$ and (ii) $\|S_t
\varphi\|_p \rightarrow 0$ for every $\varphi \in {\cal S}$, as $t \rightarrow 0$.
Let us assume these results for a moment and complete the proof.

Let $\epsilon > 0$ be given. By (i), there is a constant $C > 0$ such that
\begin{equation*}
\sup\limits_{0 \leq t \leq T} \|S_t f\|_p \leq C~\|f\|_p, f \in {\cal S}_p.
\end{equation*}
Choose $\varphi \in {\cal S}$, so that $\|f -\varphi\|_p \leq
(\frac{\epsilon}{2C}) $. Then,
\begin{align*}
\|S_t f\|_p &\leq \|S_t (f-\varphi)\|_p + \|S_t \varphi\|_p \\[.2pc]
&\leq \epsilon/2 + \|S_t \varphi\|_p.
\end{align*}
Now choose $\delta > 0$ such that $\|S_t \varphi\|_p \leq \epsilon/2
~\mbox{for all}~ 0 \leq t < \delta$, to get $\|S_t f\|_p < \epsilon$ for
all $0 \leq t < \delta$.

Since $S_t = {\cal F} (T_t - I) {\cal F}^{-1}$, (i) follows from the fact
that $T_t: {\cal S}_p \rightarrow {\cal S}_p$ are uniformly bounded
(Theorem~2.4a) and ${\cal F}: {\cal S}_p \rightarrow {\cal S}_p$ is a
unitary operator. The proof of (ii) is by a direct calculation when
$p=m$ is a non-negative integer.
\begin{equation*}
\|S_t \varphi\|_m = \|H^m S_t \varphi\|_0
\leq C_1 \sum\limits_{|\alpha| +|\beta| \leq 2m} \|x^\alpha \partial^\beta
S_t \varphi\|_0.
\end{equation*}
Since $S_t \varphi (x) = ({\rm e}^{-(t/2) |x|^2} -1) \varphi (x)$, by Leibniz rule
\begin{equation*}
\|x^\alpha \partial^\beta S_t \varphi\|_0 \leq \sum\limits_
{|\mu|+|\gamma| =|\beta|} C_{\mu \gamma} \|x^\alpha \partial^\mu ({\rm e}^{-(t/2) |x|^2} -1) \partial^
\gamma \varphi\|_0.
\end{equation*}
When $\mu \neq 0$, we have
\begin{equation*}
\|x^\alpha \partial^\mu ({\rm e}^{-(t/2) |x|^2} -1) \partial^\gamma \varphi\|_0
\leq C_2 t^{|\mu|} \|\varphi\|_m
\end{equation*}
and when $\mu=0$, using the elementary inequality $|1 - {\rm e}^{-u}| \leq C_3 u,
u > 0$ we get
\begin{equation*}
\|x^\alpha ({\rm e}^{-(t/2) |x|^2} -1) \partial^\gamma \varphi\|_0 \leq C_4 t
\|\varphi\|_{m+1}.
\end{equation*}
Therefore, $\|S_t \varphi\|_m \leq C~t \|\varphi\|_{m+1}$ for some
constant $C$, which shows that $\|S_t \varphi\|_m \rightarrow 0$ as $t
\rightarrow 0$. If $p$ is real and $m$ is a non-negative integer such that
$p \leq m$, we have
\begin{equation*}
\|S_t \varphi\|_p \leq \|S_t \varphi\|_m \leq Ct \|\varphi\|_{m+1}
\end{equation*}
and so $\|S_t \varphi\|_p \rightarrow 0$ as $t \rightarrow 0$ in this case
as well.\hfill \ab
\end{proof}

\end{document}